\documentclass[12pt]{article}

\usepackage[lmargin=2cm,rmargin=2cm,tmargin=3cm,bmargin=2cm]{geometry}

\usepackage{amsmath}    
\usepackage{amssymb}    
\usepackage{authblk}    
\usepackage{listings}   
\usepackage{hyperref}   
\usepackage{mathtools}  
\usepackage{xcolor}     

\allowdisplaybreaks{}

{}
{}
\let\leq=\leqslant{}
\let\geq=\geqslant{}

\renewcommand{\vec}[1]{\boldsymbol{#1}}

\begin{document}

\title{\bf Robust Optimal Contribution Selection}
\author[1,3]{Josh Fogg}
\author[2]{Jaime Ortiz}
\author[2]{Ivan Pocrni\'c}
\author[1]{J.~A.~Julian Hall}
\author[2]{Gregor Gorjanc}
\affil[1]{\small The Maxwell Institute, School of Mathematics, The University of Edinburgh, Edinburgh, UK}
\affil[2]{\small The Roslin Institute, Royal (Dick) School of Veterinary Science, The University of Edinburgh, Midlothian, UK}
\affil[3]{Corresponding author: \href{mailto:j.fogg@ed.ac.uk}{j.fogg@ed.ac.uk}}
\date{2024-12-03}


\maketitle

\abstract{\noindent Optimal contribution selection (OCS) is a selective breeding method that manages the conversion of genetic variation into genetic gain to facilitate short-term competitiveness and long-term sustainability in breeding programmes.
Traditional approaches to OCS do not account for uncertainty in input data, which is always present and challenges optimization and practical decision making.
Here we use concepts from robust optimization to derive a robust OCS problem and develop two ways to solve the problem using either conic optimization or sequential quadratic programming.
We have developed the open-source Python package \href{https://github.com/Foggalong/RobustOCS}{\texttt{robustocs}} that leverages the Gurobi and HiGHS solvers to carry out these methods.
Our testing shows favourable performance when solving the robust OCS problem using sequential quadratic programming and the HiGHS solver. \smallskip}


\section{Introduction}

In selective breeding programmes it is important to manage the conversion of genetic variation into genetic gain with care.
The efficiency of this conversion depends on the optimal use of genetically superior individuals to maximise short- and long-term response to selection~\cite{woolliams_genetic_2015}.

Given a cohort of \(n\) selection candidates, we approach this problem by first obtaining a breeding value \(\mu_i\) for each candidate, stored in an \(n\)-dimensional vector \(\vec{\mu}\), which describes the genetic value of individuals as parents of the next generation.
In addition, we obtain a genetic relationship matrix \(\vec{\Sigma}\), with coefficient \(\sigma_{ij}\) between candidates \(i\) and \(j\), which describes genetic similarity between candidates~\cite{henderson_simple_1976}.
Then we use \(\vec{\mu}\) and \(\vec{\Sigma}\) in optimal contribution selection (OCS)~\cite{woolliams_genetic_2015} to find a vector \(\vec{w}\) of the selection candidates' contributions to the next generation, with \(w_i\) denoting the contribution of the \(i^{.\text{th}}\) candidate.
Since \(\vec{w}\) is a vector of proportions its entries must sum to one, {\it i.e.\/} \(\sum_i w_i = 1\).
If in addition our cohort is split into male candidates (sires, \(\mathcal{S}\)) and female candidates (dams, \(\mathcal{D}\)), with each required to make half the contribution, then the sum-to-one constraint is superseded by two constraints \(\sum_{i\in\mathcal{S}}w_i = \frac{1}{2}\) and \(\sum_{i\in\mathcal{S}}w_i = \frac{1}{2}\).
We may also set upper bounds \(\vec{u}\) and lower bounds \(\vec{l}\) on contributions.

In OCS we maximize response to selection (by assigning higher contributions to candidates with higher breeding values) while minimizing risks due to inbreeding (by assigning lower contributions to closely related candidates)~\cite{woolliams_genetic_2015}.
This can be formulated as a multi-objective optimization problem,
\begin{equation}\label{eq:ocs-dual-objective}
    \begin{gathered}
        \min_{\vec{l}\leq\vec{w}\leq\vec{u}} \left( \frac{1}{2}\vec{w}^{T}\!\vec{\Sigma}\vec{w} \right),\ 
        \max_{\vec{l}\leq\vec{w}\leq\vec{u}} \left( \vec{w}^{T}\!\vec{\mu} \right)
        \ \text{subject to}\
        \sum_{i\in\mathcal{S}}w_i = \frac{1}{2}\ \text{and}\
        \sum_{i\in\mathcal{D}}w_i = \frac{1}{2}.
    \end{gathered}
\end{equation}
Since this is a multi-objective problem, we typically reframe it as a single objective maximization problem,
\begin{equation}\label{eq:ocs-single-objective}
    \max_{\vec{l}\leq\vec{w}\leq\vec{u}} \left( \vec{w}^{T}\!\vec{\mu} - \frac{\lambda}{2}\vec{w}^{T}\!\vec{\Sigma}\vec{w} \right)\ \text{s.t.}\ \sum_{i\in\mathcal{S}}w_i = \frac{1}{2},\ \sum_{i\in\mathcal{D}}w_i = \frac{1}{2},
\end{equation}
with a parameter \(\lambda\geq0\) that balances emphasis between the response to selection and risks from inbreeding.

We can solve this using quadratic programming (QP) for varying values of \(\lambda^{(k)}\) to find a set of Pareto optimal solutions \(\vec{w}^{(k)}\).
These are the points where, for the risk corresponding to \(\vec{w}^{(k)}\), there is no higher genetic merit, and for the genetic merit corresponding to \(\vec{w}^{(k)}\) there is no lower risk.
Higher \(\lambda^{(k)}\) values will correspond to a lower tolerance for risks from inbreeding.

It is helpful rearrange the problem into `standard form',
\begin{equation*}
    \min_{\vec{l}\leq\vec{x}\leq\vec{u}} \left(\vec{x}^{T}\!\vec{A}\vec{x} + \vec{q}^{T}\!\vec{x}\right)\ \text{subject to}\ \vec{G}\vec{x}\leq\vec{h}, \vec{M}\vec{x} = \vec{m},
\end{equation*}
with constraints governed by matrices \(\vec{G}\) and \(\vec{M}\) and associated vectors \(\vec{h}\) and \(\vec{m}\).
In doing so, note that the two vector constraints \(\sum_{i\in\mathcal{S}}w_i = \frac{1}{2}\) and \(\sum_{i\in\mathcal{D}}w_i = \frac{1}{2}\) can be represented as a single matrix constraint \(\vec{M}\vec{w} = \vec{m}\), where \(\vec{M}\) is a matrix with all entries either 0 and 1, and \(\vec{m}\) is a vector where all entries are \(\frac{1}{2}\).
The OCS problem does not involve the \(\vec{G}\vec{x}\leq\vec{h}\) constraint, though it would be straightforward to adapt the theory to include this.
Hence, OCS in standard form is then
\begin{equation}\label{eq:ocs-standard-form}
    \max_{\vec{l}\leq\vec{w}\leq\vec{u}} \left( \vec{w}^{T}\!\vec{\mu} - \frac{\lambda}{2}\vec{w}^{T}\!\vec{\Sigma}\vec{w} \right)\ 
    \text{subject to}\ \vec{M}\vec{w} = \vec{m}.
\end{equation}
Similarly, the theory presented may easily be adapted for use with single-sex cohorts or to add general linear constraints, \(\vec{A}\vec{w} = \vec{b}\), where \(\vec{A}\) is some matrix and \(\vec{b}\) some vector.


\section{Incorporating uncertainty}\label{sect:AssumptionOfOCS}

Breeding values for a given trait are estimated using information on phenotypes and genetic relationships among individuals in a population~\cite{mrode_genetic_2014}.
This estimation involves the calculation of conditional expectation \(\bar{\vec{\mu}} = \mathbb{E}\left(\vec{\mu} | \vec{y}\right)\) (referred to as estimated breeding value, EBV) and conditional covariance \(\vec{\Omega} = \operatorname{Var}\left(\vec{\mu} | \vec{y}\right)\), where \(\vec{y}\) is a vector of observed phenotypes, and expectation and covariance are with respect to an assumed estimation model.
These conditional quantities can be obtained using closed-form, iterative, or sampling approaches~\cite{mrode_genetic_2014,hickey_estimation_2009,garcia-cortes_estimation_1995,fouilloux_sampling_2001}.
The conditional variance of EBV decreases as more information becomes available over time.
At the time of selection, the amount of information is usually limited, and hence there is uncertainty associated with the EBV.
In addition, EBV are also correlated due to conditioning on observed phenotypes.
Thus any solution to~\eqref{eq:ocs-dual-objective} or~\eqref{eq:ocs-single-objective} is contingent on the implicit and convenient assumption that the values of \(\vec{\mu}\) and \(\vec{\Sigma}\) are known exactly; in particular, we implicitly assume that EBV \((\bar{\vec{\mu}})\) are equal to true breeding values \((\vec{\mu})\).
Unfortunately, this is not the case, and there is a certain probability that they will differ.
Hence, we should account for the uncertainty and correlation of EBV in OCS.

We use robust optimization to address this problem.
To this end, we reformulate~\eqref{eq:ocs-standard-form} as a bilevel optimization problem
\begin{equation*}
    \max_{\vec{l}\leq\vec{w}\leq\vec{u}} \left( \min_{\vec{\mu}\in U_{\vec{\mu}}} \left( \vec{w}^{T}\!\vec{\mu} \right) - \frac{\lambda}{2}\vec{w}^{T}\!\vec{\Sigma}\vec{w} \right)\ 
    \text{subject to}\ M\vec{w} = \vec{m},
\end{equation*}
where \(U_{\vec{\mu}}\) is some uncertainty set of values for
\(\vec{\mu}\) that we deem allowable.
In bilevel optimization this is called a follower-leader setup.
We solve the inner problem first, with its solution impacting the outer problem.

While bilevel optimization is a wide area of research, in our setup we can solve the inner problem explicitly based on our choice of \(U_{\vec{\mu}}\)~\cite{yin_practical_2021}.
For practical reasons relating to continuity and differentiability, we define a `quadratic' uncertainty set\footnotemark{}:
\begin{equation*}
    U_{\vec{\mu}} := \left\lbrace \vec{\mu} :\ {(\vec{\mu}-\bar{\vec{\mu}})}^{T}\!\vec{\Omega}^{-1}\!(\vec{\mu}-\bar{\vec{\mu}}) \leq \kappa^2 \right\rbrace,
\end{equation*}
which bounds the uncertainty in a ball about \(\bar{\vec{\mu}}\), with parameter \(\kappa\) controlling the tolerance for uncertainty.
\footnotetext{~
  We could instead use \(U_{\vec{\mu}} := \left\lbrace \vec{\mu} :\ | \mu_i - \bar{\mu}_i| \leq\xi_i,\ \forall i \right\rbrace\) (a `box' uncertainty set), where \(\vec{\xi}\) is some upper bound~\cite{heckel_insights_2016}. However, this does not utilise the known posterior covariance \(\vec{\Omega}\).
}
Using the quadratic uncertainty set means that inner problem \(\min_{\vec{\mu}\in U_{\vec{\mu}}} (\vec{w}^{T}\!\vec{\mu})\) in standard form becomes
\begin{equation}\label{eq:rocs-inner-problem}
    \min_{\vec{\mu}} \left(\vec{w}^{T}\!\vec{\mu}\right) \quad\text{subject to}\quad \kappa^2 - {(\vec{\mu}-\bar{\vec{\mu}})}^{T}\!\vec{\Omega}^{-1}\!(\vec{\mu}-\bar{\vec{\mu}}) \geq 0.
\end{equation}

\subsection{Intuitive example}\label{sec:example-intuition}

To get an idea of why this works intuitively, we present a toy problem. Consider a cohort with three candidates represented by the following data:
\begin{equation*}
    \bar{\vec{\mu}} = \begin{bmatrix}
        1 \\
        2 \\
        1
    \end{bmatrix},
    \qquad \vec{\Omega} = \begin{bmatrix}
        1/9 & 0 & 0 \\
        0   & 4 & 0 \\
        0   & 0 & 1
    \end{bmatrix},
    \qquad \vec{\Sigma} = \begin{bmatrix}
        1 & 0 & 0 \\
        0 & 1 & 0 \\
        0 & 0 & 1
    \end{bmatrix},
    \qquad \mathcal{S} = \lbrace 1, 2 \rbrace,
    \qquad \mathcal{D} = \lbrace 3 \rbrace.
\end{equation*}

By inspection we can see that since we only have one dam, it will make up the full 50\% that dams contribute, while the 50\% that sires contribute is split between candidates \#1 and \#2. The standard solution which does not take into account uncertainty in breeding values is \(\vec{w} = {\begin{bmatrix} 0 & 0.5 & 0.5 \end{bmatrix}}^{T}\) since the second candidate has the better mean expected breeding value.

However, although the mean expected breeding value for the first candidate is smaller, its variance is far smaller than that of the second. Hence, although the mean expected breeding value for the second candidate is higher, its true expected breeding value may be far lower than that of the first.

Consider the constraint from the inner problem~\eqref{eq:rocs-inner-problem},
\begin{equation*}
    \kappa^2 \geq {\left( \begin{bmatrix}
        \mu_1 \\ \mu_2 \\ \mu_3
    \end{bmatrix} - \begin{bmatrix}
        1 \\ 2 \\ 1
    \end{bmatrix}\right)}^{T}\!\begin{bmatrix}
        9 &   0 & 0 \\
        0 & 1/4 & 0 \\
        0 &   0 & 1
    \end{bmatrix}{\left( \begin{bmatrix}
        \mu_1 \\ \mu_2 \\ \mu_3
    \end{bmatrix} - \begin{bmatrix}
        1 \\ 2 \\ 1
    \end{bmatrix}\right)}
    = 9{(\mu_1 - 1)}^2 + \frac{1}{4}{(\mu_2 - 2)}^2 + {(\mu_3 - 1)}^2
\end{equation*}
This show that, for varying values of \(\kappa\), the boundaries of the uncertainty sets are ellipsoids centred at \(\begin{bmatrix} 1 & 2 & 1 \end{bmatrix}^{T}\). For \(\kappa = 1\) for example, by examining the equation or plotting the ellipsoid, we can see that
\begin{equation*}
    \frac{2}{3} \leq \mu_1 \leq \frac{4}{3},\qquad 0 \leq \mu_2 \leq 4,\quad\text{ and }\quad 0 \leq \mu_3 \leq 2.
\end{equation*}
This gives a `worst case' expected breeding value of \(\frac{2}{3}\) for the first candidate and 0 for the second and third candidates. Note that the use of a quadratic uncertainty set (a `ball'), rather than a box uncertainty set, means that all worst cases cannot be achieved simultaneously. This demonstrates how robust optimization encourages pessimism, but not excessively.

When solving the robust optimization problem, for small values of \(\kappa\),\footnotemark{} the solution which didn't take into account uncertainty (\({\begin{bmatrix} 0 & 0.5 & 0.5 \end{bmatrix}}^{T}\)) is obtained: the uncertainty set is too small for the `worst case' expecting breeding value of the second candidate to be less than the mean expected breeding value of the first. However, as \(\kappa\) increases, the robust solution shifts towards \({\begin{bmatrix} 0.5 & 0 & 0.5 \end{bmatrix}}^{T}\), {\it i.e.\/} a solution which puts a higher weight on the `safer' sire candidate.

\footnotetext{~Observe that if \(\kappa = 0\), the robust problem full reduces to the standard non-robust problem.}

\subsection{General solution}\label{sect:GeneralSolution}

Returning to bilevel optimization problem's inner problem~\eqref{eq:rocs-inner-problem}, we know this is convex because \(\vec{\Omega}\) is a positive definite matrix and \(\vec{w}^{T}\!\vec{\mu}\) is a linear function.
With these conditions met, the first-order necessary conditions (Karush-Kuhn-Tucker) are necessary and sufficient for optimality.
If we define a Lagrange multiplier \(\rho\in\mathbb{R}\), the conditions for this problem are:
\begin{align}
    \nabla_{\vec{\mu}}L(\vec{\mu}, \rho) = 0 \quad&\Rightarrow\quad \nabla_{\vec{\mu}}\left( \vec{w}^{T}\!\vec{\mu} - \rho\big( \kappa^2 - {(\vec{\mu}-\bar{\vec{\mu}})}^{T}\!\vec{\Omega}^{-1}\!(\vec{\mu}-\bar{\vec{\mu}}) \big) \right) = 0,\label{eq:RobKKT1} \\
    c(\vec{\mu}) \geq 0 \quad&\Rightarrow\quad \kappa^2 - {(\vec{\mu}-\bar{\vec{\mu}})}^{T}\!\vec{\Omega}^{-1}\!(\vec{\mu}-\bar{\vec{\mu}}) \geq 0, \nonumber\\ 
    \rho \geq 0 \quad&\Rightarrow\quad \rho\geq0, \nonumber\\
    \rho \cdot c(\vec{\mu}) = 0 \quad&\Rightarrow\quad \rho\left(\kappa^2 - {(\vec{\mu}-\bar{\vec{\mu}})}^{T}\!\vec{\Omega}^{-1}\!\vec{\mu}-\bar{\vec{\mu}}) \right) = 0. \label{eq:RobKKT4}
\end{align}
From~\eqref{eq:RobKKT1} we have that:
\begin{equation}\label{eq:RobKKT5}
    \vec{w} + \rho 2\vec{\Omega}^{-1}\!(\vec{\mu}-\bar{\vec{\mu}}) = 0 \quad\Rightarrow\quad \vec{\mu} - \bar{\vec{\mu}} = -\frac{1}{2\rho}\vec{\Omega}\vec{w},
\end{equation}
which when substituted into~\eqref{eq:RobKKT4} gives:
\begin{align*}
    &\phantom{\Rightarrow}\quad\ \ \rho\left( \kappa^2 - \big(-\frac{1}{2\rho}\vec{\Omega} \vec{w}\big)^{T}\!\vec{\Omega}^{-1}\!\big(-\frac{1}{2\rho}\vec{\Omega} \vec{w}\big) \right) = 0 \\
    &\Rightarrow\quad \rho\kappa^2 - \frac{1}{4\rho}\vec{w}^{T}\!\vec{\Omega}^{T}\!\vec{\Omega}^{-1}\!\vec{\Omega}\vec{w} = 0 \\
    &\Rightarrow\quad \rho^2\kappa^2 = \frac{1}{4}\vec{w}^{T}\!\vec{\Omega} \vec{w}\quad \text{(since \(\vec{\Omega}\) is symmetric)} \\
    &\Rightarrow\quad \rho\kappa = \frac{1}{2}\sqrt{\vec{w}^{T}\!\vec{\Omega}\vec{w}}\quad \text{(since \(\rho,\kappa\geq0\))} \\
    &\Rightarrow\quad \frac{1}{2\rho} = \frac{\kappa}{\sqrt{\vec{w}^{T}\!\vec{\Omega}\vec{w}}}.
\end{align*}
Substituting this back into~\eqref{eq:RobKKT5}, we find that
\begin{equation*}
    \vec{\mu} - \bar{\vec{\mu}} = -\frac{\kappa}{\sqrt{\vec{w}^{T}\!\vec{\Omega}\vec{w}}}\vec{\Omega}\vec{w} \quad\Rightarrow\quad \vec{\mu} = \bar{\vec{\mu}} - \frac{\kappa}{\sqrt{\vec{w}^{T}\!\vec{\Omega}\vec{w}}}\vec{\Omega}\vec{w}
\end{equation*}
is the solution for the inner problem.
After substituting this back into the outer problem and rationalising, we obtain
\begin{equation}\label{eq:robust-ocs-sqrt}
    \max_{\vec{l}\leq\vec{w}\leq\vec{u}} \left(\vec{w}^{T}\bar{\vec{\mu}} - \kappa\sqrt{\vec{w}^{T}\vec{\Omega}\vec{w}} - \frac{\lambda}{2}\vec{w}^{T}\vec{\Sigma} \vec{w}\right)\ \text{s.t.}\ \vec{M}\vec{w} = \vec{m},
\end{equation}
where \(\kappa\in\mathbb{R}\) is the robust optimization parameter controlling the tolerance for uncertainty.
Since our objective has gained a square root term, \eqref{eq:robust-ocs-sqrt} is no longer a quadratic problem.

\subsection{Example solution}

We return to our example cohort modelled by variables 
\begin{equation*}
    \bar{\vec{\mu}} = \begin{bmatrix}
        1 \\
        2 \\
        1
    \end{bmatrix},
    \qquad \vec{\Omega} = \begin{bmatrix}
        1/9 & 0 & 0 \\
        0   & 4 & 0 \\
        0   & 0 & 1
    \end{bmatrix},
    \qquad \vec{\Sigma} = \begin{bmatrix}
        1 & 0 & 0 \\
        0 & 1 & 0 \\
        0 & 0 & 1
    \end{bmatrix},
    \qquad \mathcal{S} = \lbrace 1, 2 \rbrace,
    \qquad \mathcal{D} = \lbrace 3 \rbrace,
\end{equation*}
which we wish to find the we solve~\eqref{eq:robust-ocs-sqrt} for \(\lambda = 0.1\) and \(\kappa = 1\).

Since we only have one dam, candidate three, it follows in any solution \(w_3 = 0.5\). Similarly, since there are only two sires and their contributions must sum to half, it follows that \(w_2 = 0.5 - w_1\). Thus we can restate~\eqref{eq:robust-ocs-sqrt} in this case as
\begin{multline*}
    \max_{0 \leq w_1 \leq 0.5} \left(\begin{bmatrix}
        w_1 \\ 0.5 - w_1 \\ 0.5
    \end{bmatrix}^{T}\!\begin{bmatrix}
        1 \\ 2 \\ 1
    \end{bmatrix} - \sqrt{\begin{bmatrix}
        w_1 \\ 0.5 - w_1 \\ 0.5
    \end{bmatrix}^{T}\!\begin{bmatrix}
        1/9 & 0 & 0 \\
        0   & 4 & 0 \\
        0   & 0 & 1
    \end{bmatrix}\begin{bmatrix}
        w_1 \\ 0.5 - w_1 \\ 0.5
    \end{bmatrix}}\right. \\ \left. - \frac{0.1}{2}\begin{bmatrix}
        w_1 \\ 0.5 - w_1 \\ 0.5
    \end{bmatrix}^{T}\!\begin{bmatrix}
        w_1 \\ 0.5 - w_1 \\ 0.5
    \end{bmatrix}\right),
\end{multline*}
which after a lot of rearranging gives us
\begin{equation*}
    \max_{0 \leq w_1 \leq 1} \left( -0.1w_1^2 - 0.95w_1^{\phantom{2}} + 1.475 - \sqrt{\frac{37}{9}w_1^2 - 4w_1^{\phantom{2}} + 1.25w_1^{\phantom{2}}} \right).
\end{equation*}
This is a differentiable equation in a single variable, so we can employ basic calculus to find the maximum value. We find that the robust solution is \(\vec{w} = \begin{bmatrix} 0.3359 & 0.1641 & 0.5 \end{bmatrix}^{T}\).

To see how this aligns with our intuition, observe that:
\begin{itemize}
    \item The solution ignoring uncertainty is \(\vec{w} = {\begin{bmatrix} 0 & 0.5 & 0.5 \end{bmatrix}}^{T}\), with objective value \(1.475\), coming from \(\vec{w}^{T}\!\bar{\vec{\mu}} = 1.5\) and \(-\frac{\lambda}{2}\vec{w}^{T}\!\Sigma\vec{w} = -0.025\).
    \item Taking this solution and substituting it into the inner problem, we see that~\eqref{eq:rocs-inner-problem} becomes
    \begin{equation*}
        \min_{\vec{\mu}} \left(\frac{\mu_2}{2} + \frac{\mu_3}{2}\right) \ \text{ subject to }\ 9{(\mu_1 - 1)}^2 + \frac{1}{4}{(\mu_2 - 2)}^2 + {(\mu_3 - 1)}^2 \leq 1.
    \end{equation*}
    Since \(\mu_1\) is not in the objective, we may set its value as high as possible to meet the constraint, ultimately giving \(\vec{\mu} = {\begin{bmatrix} \frac{4}{3} & 0 & 0 \end{bmatrix}}^{T}\). With this though, since \(\vec{w}^{T}\!\vec{\mu} = 0\), the objective value of the outer problem now drops to \(-0.025\).
    
    \item In fact the robust solution is \(\vec{w} = \begin{bmatrix} 0.3359 & 0.1641 & 0.5 \end{bmatrix}^{T}\) (found by directly solving or using one the the methods discussed below), corresponding to an objective value of \(0.5361\) and \(\vec{\mu} = \begin{bmatrix} 0.9387 & 0.9215 & 0.1782 \end{bmatrix}^{T}\).
\end{itemize}
Thus, although the robust solution objective value of 0.5361 is well down on the non-robust objective of 1.475, it's rather better than the `worst case' of -0.025.

We were only able to find the robust solution directly here because it was a toy problem. To tackle any practical problem we will need other methods.

\section{Direct conic optimization}\label{sect:ConicOptimization}

Working with this new form~\eqref{eq:robust-ocs-sqrt} is not ideal, and short of using a general non-linear solver it is not a problem most off-the-shelf software can handle.
However, we can make it tractable by adapting it into the form of a conic optimization problem.
If we define a real auxiliary variable \(z \geq 0\) such that \(z \geq \sqrt{\vec{w}^{T}\!\vec{\Omega}\vec{w}}\), then the problem becomes:
\begin{equation*}
    \max_{\substack{\vec{l}\leq\vec{w}\leq\vec{u}, \\ z\geq0}} \left(\vec{w}^{T}\!\bar{\vec{\mu}} - \kappa z - \frac{\lambda}{2}\vec{w}^{T}\!\vec{\Sigma}\vec{w}\right)\ \text{s.t.}\ z \geq \sqrt{\vec{w}^{T}\!\vec{\Omega}\vec{w}},\ \vec{M}\vec{w} = \vec{m}.
\end{equation*}
Since \(z > 0\), \(\kappa \geq 0\), and we are maximizing an objective containing
`\(-\kappa z\)', this term of the objective will be biggest when \(z\) is smallest.
This happens precisely when it attains its lower bound from the constraint
\(z \geq \sqrt{\vec{w}^{T}\!\vec{\Omega}\vec{w}}\),
hence \(z\) will push downwards and our relaxation is justified.

However, the presence of the square root still proves problematic for most optimization software, so we further make use of both \(z\) and \(\sqrt{\vec{w}^{T}\!\vec{\Omega}\vec{w}}\) being positive to note that \(z \geq \sqrt{\vec{w}^{T}\!\vec{\Omega}\vec{w}}\) can be squared on both sides:
\begin{equation}\label{eq:robust-ocs-conic}
    \max_{\substack{\vec{l}\leq \vec{w}\leq \vec{u}, \\ z\geq0}} \left(\vec{w}^{T}\!\bar{\vec{\mu}} - \kappa z - \frac{\lambda}{2}\vec{w}^{T}\!\vec{\Sigma}\vec{w} \right)\ \text{s.t.}\ z^2 \geq \vec{w}^{T}\!\vec{\Omega}\vec{w},\ M\vec{w} = \vec{m}.
\end{equation}
Here \(z^2 \geq \vec{w}^{T}\vec{\Omega} \vec{w}\) represents a cone constraint, hence we are solving a second order conic optimization problem.
This opens up the possibility of using standard optimizaton software.

For example, provided we have the variables loaded into Python, we can use Gurobi's modelling language~\cite{gurobi2024} to solve such a problem, as indicated with the code below. 

\begin{minipage}{\hsize}%
    \begin{lstlisting}
model = gp.Model("robust-genetics")  # initializes

w = model.addMVar(shape=dimension, name="w", lb=lower_bound, 
                  ub=upper_bound, vtype=gp.GRB.CONTINUOUS)
z = model.addVar(lb=0.0, name="z")

# NOTE lambda is `lam' to avoid Python conflicts
model.setObjective(
    w.T@mubar - (lam/2)*w.T@(sigma@w) - kappa*z,
    gp.GRB.MAXIMIZE
)

model.addConstr(M@w == m, name="sum-to-half")
model.addConstr(z**2 >= w.T@omega@w, name="uncertainty")

model.optimize()
    \end{lstlisting}
\end{minipage}

\section{Sequential quadratic programming}\label{sect:SQP}

For our other approach, from~\eqref{eq:robust-ocs-conic} we note \(f(\vec{w}) = \sqrt{\vec{w}^{T}\!\vec{\Omega}\vec{w}}\) is a cone with a derivative discontinuity at the origin.
For \(\vec{w} \neq \vec{0}\) (which we can guarantee since \(w_i \geq 0\) and \(\sum_i w_i = 1\)) we have:
\begin{equation*}
    \nabla f(\vec{w}) = \frac{1}{\sqrt{\vec{w}^{T}\!\vec{\Omega}\vec{w}}}\vec{\Omega}\vec{w}.
\end{equation*}
This means we can approximate the cone using a series of tangent planes (with the \(k^{\text{th}}\) denoted \(P_k\)) of the form
\begin{equation*}
    z \geq \frac{1}{\sqrt{\vec{w}^{T}_{(k)}\!\vec{\Omega}\vec{w}^{\phantom{T}}_{(k)}}} {\left(\vec{\Omega}\vec{w}_{(k)}\right)}^{T} \vec{w},
\end{equation*}
for a set of \(K\) points \(\vec{w}_{(0)}, \vec{w}_{(1)}, \ldots, \vec{w}_{(K)}\). To find \(\vec{w}_{(K)}\) we use sequential quadratic programming (SQP), which solves the quadratic problem:
\begin{equation}\label{eq:robust-ocs-sqp}
    \begin{gathered}
        \max_{\substack{\vec{l}\leq\vec{w}\leq\vec{u}, \\ z\geq0}} \left(\vec{w}^{T}\!\bar{\vec{\mu}} - \kappa z - \frac{\lambda}{2}\vec{w}^{T}\!\vec{\Sigma}\vec{w} \right)\ \text{s.t.}\ M\vec{w} = \vec{m},\ \\
        z \geq \frac{1}{\sqrt{\vec{w}^{T}_{(k)}\!\vec{\Omega}\vec{w}^{\phantom{T}}_{(k)}}} {\left(\vec{\Omega}\vec{w}_{(k)}\right)}^{T} \vec{w}\ \text{ for }\ k = 0, \ldots, (K-1).
    \end{gathered}
\end{equation}

This is a regular (convex) QP problem and, as before, we can use Gurobi's modelling language using code such as that below to solve this.
An important advantage of SQP is that it only relies on QP, so a greater number of solvers may be used. 

\begin{minipage}{\hsize}%
    \begin{lstlisting}
model = gp.Model("robust-genetics-sqp")  # initializes

w = model.addMVar(shape=dimension, name="w", lb=lower_bound, 
                  ub=upper_bound, vtype=gp.GRB.CONTINUOUS)
z = model.addVar(lb=0.0, name="z")

# NOTE lambda is `lam' to avoid Python conflicts
model.setObjective(
    w.T@mubar - (lam/2)*w.T@(sigma@w) - kappa*z,
    gp.GRB.MAXIMIZE
)

model.addConstr(M@w == m, name="sum-to-half")

for i in range(max_iterations):
    model.optimize()

    # z coefficient for the new constraint
    w_star = np.array(w.X)
    alpha = sqrt(w_star.T@omega@w_star)

    # if gap between z and w'Omega w converges, done
    if abs(z.X - alpha) < robust_gap_tol: break

    # add a new plane to the cone's approximation
    model.addConstr(alpha*z >= w_star.T@omega@w, name=f"P{i}")
    \end{lstlisting}
\end{minipage}

\section{Implementation in Python}\label{sect:Python}

Though freely available to academics, Gurobi is commercial optimization software with license fees which are prohibitive to many stakeholders.
To democratise this work, we provide a solution that leverages open-source software.
One such tool is \href{https://highs.dev/}{HiGHS}~\cite{huangfu_parallelizing_2018}, released under the \href{https://choosealicense.com/licenses/mit}{MIT license}.
It can solve convex QP problems so is applicable to robust OCS
using SQP~\eqref{eq:robust-ocs-sqp}.

For ease of use and testing, we implemented these methods as the Python package \texttt{robustocs}, available via \href{https://pypi.org/project/robustocs/}{PyPI}.
It accesses the Gurobi API using \texttt{\href{https://pypi.org/project/gurobipy}{gurobipy}} and the HiGHS API using \texttt{\href{https://pypi.org/project/highspy}{highspy}}, also using \texttt{\href{https://numpy.org/}{numpy}} for standard linear algebra tools and \texttt{\href{https://scipy.org/}{scipy}} for handling sparse matrix objects.
The \texttt{robustocs} package is released under the MIT license, with development and documentation at \href{https://github.com/Foggalong/robustocs}{github.com/Foggalong/robustocs}.

Suppose we have a breeding cohort with genetic relationship matrix \(\vec{\Sigma}\) and EBV vector \(\bar{\vec{\mu}}\) with associated covariance matrix \(\vec{\Omega}\).
If these are saved in files of appropriate formats, solving the robust OCS using \texttt{robustocs} can be done by running the code below in Python. 
\texttt{robustocs} can also take inputs via NumPy or SciPy objects and has more granular functions to control over the solver, method, and associated parameters ({\it e.g.\/} maximum solve time). 

\begin{minipage}{\hsize}%
    \begin{lstlisting}
import robustocs as rocs

selection, objective, merit, coancestry = rocs.solveROCS(
    sigma_filename="cohort-relationships.txt",
    mu_filename="breeding-means.txt",
    omega_filename="breeding-variances.txt",
    sex_filename="cohort-sexes.txt",
    method='robust', lam=0.5, kappa=1
)
    \end{lstlisting}
\end{minipage}

To evaluate how the methods performed when solving robust OCS, we simulated an example dataset using AlphaSimR~\cite{gaynor_alphasimr_2021}.
The simulation mimicked a breeding programme over 10 generations with 1000 individuals in each.
In each generation, we selected the best 25 males (out of 500) as fathers based on their phenotypes and mated them with all 500 females from the previous generation and all 500 females from the current generation.
These matings produced 1000 selection candidates for the next generation, in total 10,000 across 10 generations.
We fitted the pedigree-based linear mixed model~\cite{mrode_genetic_2014} to the simulated data and the known \(\vec{\Sigma}\) from pedigree.
From this model fit, we obtained 1000 samples from the posterior distribution \(p\left(\vec{\mu} | \vec{y}\right)\) to estimate \(\bar{\vec{\mu}}\) and \(\vec{\Omega}\).
While \(\vec{\Omega}\) is also estimated, modelling its uncertainty is beyond the scope of this study and will likely have diminishing returns.

Using this, we compared how the methods performed across Gurobi and HiGHS as the size of the cohort increased.\footnotemark{}
The results in Table~\ref{tab:large-examples} show favourable execution time for using HiGHS compared to Gurobi for the standard OCS problem~\eqref{eq:ocs-standard-form}.

\begin{table}[!h]
    \centering
    \begin{tabular}{rrrrrr}
        \textbf{Size} \(n\) & \textbf{Gurobi}~\eqref{eq:ocs-standard-form} & \textbf{HiGHS}~\eqref{eq:ocs-standard-form} & \textbf{Gurobi}~\eqref{eq:robust-ocs-conic} & \textbf{Gurobi}~\eqref{eq:robust-ocs-sqp} & \textbf{HiGHS}~\eqref{eq:robust-ocs-sqp} \\
            4 &  0.003 &  0.000 & 0.005 &    0.017 &   0.006 \\
           50 &  0.004 &  0.001 & 0.010 &    0.055 &   0.018 \\
         1000 &  0.676 &  0.204 & 2.750 &   26.400 &   1.680 \\
        10000 & 86.300 & 25.800 &   DNF & 1560.000 & 106.000
    \end{tabular}
    \caption{Time in seconds (to 3 s.f.) to solve standard or robust OCS problems with each method implemented with Gurobi or HiGHS against increasing problem size:~\eqref{eq:ocs-standard-form} is standard OCS,~\eqref{eq:robust-ocs-conic} is robust OCS using conic optimization, and~\eqref{eq:robust-ocs-sqp} is robust OCS using SQP}\label{tab:large-examples}
\end{table}

For the robust OCS problem solved with conic optimization~\eqref{eq:robust-ocs-conic}, Gurobi crashed without displaying an error message.
Lastly, HiGHS outperformed Gurobi by an order of magnitude when solving the robust OCS problem with SQP~\eqref{eq:robust-ocs-sqp}, and proved to be a scalable approach.

\footnotetext{~Ran on a HP~Elitebook with 15.0 GiB of memory and an Intel\textregistered~Core\textsuperscript{\tiny TM} i5-8350U CPU @ 1.70GHz \(\times\) 8, with Ubuntu 22.04.4 LTS (64-bit), \texttt{robustocs} 0.2.1, \texttt{gurobipy} 11.0.3, \texttt{highspy} 1.7.2, \texttt{numpy} 1.21, and \texttt{scipy} 1.8.0.}

\section{Conclusion}

We have proposed two robust optimization models to account for uncertainty in optimal contribution selection problems, and implemented these in a well-documented Python package \texttt{robustocs}.
The package can use HiGHS or Gurobi APIs to solve robust OCS, with HiGHS demonstrating the better performance.

\section{Acknowledgements}

The authors acknowledge support from the Harmonised Impact Acceleration Account (RobustOptimApp, H029/PV156) to JF, JAJH, and GG, the BBSRC ISP grant to The Roslin Institute (BBS/E/D/30002275, BBS/E/RL/230001A, BBS/E/RL/230001C) for IP and GG, the Roslin Institute PhD fellowship to JO, and The University of Edinburgh. For the purpose of open access, the authors have applied a CC BY public copyright license to any author-accepted manuscript version arising from this submission.


\bibliographystyle{ieeetr}
\bibliography{references}  


\end{document}